\RequirePackage{fix-cm}
\documentclass{svjour3}                     
\smartqed  
\usepackage{graphicx}
\usepackage{amsfonts,enumerate,amsmath,amssymb}
\def\qed{\strut\hfill $\Box$}
\def\para#1{\vskip .4\baselineskip\noindent{\bf #1}}
\begin{document}

\title{Stochastic averaging for non-Lipschitz multi-valued stochastic differential equations driven by $G$-Brownian motion
}

\titlerunning{Stochastic averaging for multi-valued SDEs driven by $G$-Brownian motion}        

\author{Min Han  \and Bin Pei$^*$ 
}


\institute{
	Min Han \at
              School of Mathematics and Statistics, Northwestern Polytechnical University, Xi'an, 710072, China \\
              \email{minhan2019@hotmail.com}           
                        \and
             $^*$ Corresponding author: Bin Pei \at 
             School of Mathematics and Statistics, Northwestern Polytechnical University, Xi'an, 710072, China \\
             \email{binpei@hotmail.com}
}

\date{Received: date / Accepted: date}

\maketitle

\begin{abstract}
In this paper, we prove the validity of an averaging principle for multi-valued stochastic differential equations (MSDEs) driven by $G$-Brownian motion with non-Lipschitz coefficients. The convergence theorem between the solution of the averaged MSDEs and original one was obtained in the sense of $p$-th moments and also in capicity. Finally, one example is presented to illustrate our theory.
\end{abstract}

\section{Introduction}
This paper considers the following multi-valued stochastic differential equations (MSDEs) driven by $G$-Brownian motion on $\mathbb{R}^d$:
\begin{eqnarray} \label{orig1}
dX(t)+\partial\varphi(X(t)) \ni f(t,X(t))dt+g(t,X(t))d\langle B\rangle_{t}+\sigma(t,X(t))dB_t,
\end{eqnarray}
where the initial condition $X(0)= \xi\in \overline{{\rm Dom}(\varphi)}$ with $\hat{\mathbb{E}}[{|\xi|^2}] < \infty$ and $\partial\varphi$ is the subdifferential operator associated to $\varphi$ which is a lower and semi-continuous (l.s.c) function on $\mathbb{R}^d$. $(\langle B  \rangle_t)_{t\geq 0}$ is the quadratic variation process of $G$-Brownian motion $(B_t)_{t\geq 0}$.

C\'{e}pa~\cite{Cepa1998Probleme}~firstly studied the existence and uniqueness theorem for MSDEs driven by Brownian motion. Ren, Xu and Zhang \cite{Ren2014} investigated the existence, uniqueness for MSDEs driven by continuous semimartingales. Recently, Ren, Wang and Huang \cite{Ren2017} proved the existence and uniqueness of a solution for a class of MSDEs driven by $G$-Brownian motion by means of the Yosida approximation
method. Let $\sigma=0$ and $\partial\varphi=A$ ($A$ is a multi-valued maximal montone operator) in Eq. (\ref{orig1}), Ngoran and Modeste \cite{Ngoran2001} studied the averaging principle of MSDEs driven by Brownian motion. Xu and Liu \cite{Xu2014b} removed the integrability condition about the multi-valued maximal montone operator in \cite{Ngoran2001} and obtained the convergence result between the averaged MSDE and the original one. Very recently, Guo and Pei \cite{Guo2018} established an averaging principle for MSDEs driven by Poisson
point processes. Later, Mao, Hu, You and Mao \cite{Mao2019} extended the result in \cite{Guo2018} and studied the averaging principle for MSDEs with jumps under non-Lipschitz condition.  Refer to \cite{Freidlin1998,Khasminskii1968} for more results on averaging principle.

Consider a Black-Scholes like market with uncertain volatilities, the basic securities consist of 2 assets, a riskless one, the bond, and a risky securities, the stock. Their prices are governed by
\begin{eqnarray*}
	dS^{0}_{t}=rS^{0}_{t}dt,
\end{eqnarray*}
for the bond with a constant interest rate $r>0$.
The price $S^{1}_{t}$ satisfied Black-Scholes-Merton's model driven by $G$-Brownian motion:
\begin{eqnarray*}
	dS^{1}_{t}=b_tS^{1}_{t}dt+\beta_tS^{1}_{t}d\langle B\rangle_{t}+\sigma_tS^{1}_{t}dB_t,S^{1}_{0}=x,
\end{eqnarray*}
for the stock.~$b_t,\beta_t,\sigma_t$ are assumed
to be deterministic functions of $t$. $(\langle {B} \rangle_{t})_{t\geq 0}$ is the quadratic variation process of $G$-Brownian motion $(B_t )_{t\geq 0}$.
In this market, the stock price evolution does not only
involve risk modeled by the noise part but also ambiguity about the risk due
to the unknown deviation of the process $B$ ($G$-Brownian motion) from its mean. Thus, $G$-Brownian motion is useful to be applied to measure the mean-uncertainty of risky positions. Peng (see e.g., \cite{Peng2007,Peng2008a,Peng2008b}) firstly studied the notion of sublinear expectation space and the fundamental theory of $G$-expectation and $G$-Brownian motion. Refer to \cite{Gao2009,Gao2010,lin2013} for more results about $G$-Brownian motion.

However, to the best of our knowledge, there are few literature on using the
averaging methods to obtain the approximate solutions to MSDEs in sublinear expectation space.  He et. al., \cite{He2019} proved the averaging principle for neutral functional SDEs driven by $G$-Brownian motion under the usual Lipschitz condition, unfortunately, the method proposed in \cite{He2019} failed to the MSDEs case. Moreover, most of the scholars considered the following non-Lipschitz condition (see e. g. \cite{Cao2005}):
\begin{enumerate}
	\item [(H0)] For any $t\in [0,T]$,
	\begin{eqnarray*}
	|f(t, x)-f(t, y)|^2+|\sigma(t, x)-\sigma(t, y)|^{2}\leq  \lambda(t) \gamma (\sup_{0 \leq r \leq t}|x(r)-y(r)|^{2})
\end{eqnarray*}
	where  $\lambda(t)>0$  is locally square integrable and $\gamma(\cdot): \mathbb{R}_{+} \rightarrow \mathbb{R}_{+} $  is concave nondecreasing continuous function such that $\gamma(0)=0$ and $$ \int_{0+}\frac{1}{\gamma(u)} d u=\infty.$$
\end{enumerate}

In this paper, we will study the stochastic averaging to Eq. (\ref{orig1}) under non-Lipschitz condition as follows:
\begin{enumerate}
	\item [(H1)]  The functions $f:
	\mathbb{R}_{+} \times \mathbb{R}^d \rightarrow \mathbb{R}^d$, $g : \mathbb{R}_{+}  \times
	\mathbb{R}^d \rightarrow \mathbb{R}^{d \times n}$ and $\sigma : \mathbb{R}_{+}  \times
	\mathbb{R}^d \rightarrow \mathbb{R}^{d \times n}$ are continuous and for any $x$, $y \in \mathbb{R}^d$, we have
	\begin{eqnarray*}
		&&\left(\mathbf{H}_{f}\right) \quad|f(t, x)-f(t, y)| \leq \lambda(t)|x-y| \kappa_{1}(|x-y|),\cr
		&&\left(\mathbf{H}_{g}\right) \quad|g(t, x)-g(t, y)| \leq \lambda(t)|x-y| \kappa_{2}(|x-y|),\cr
		&&\left(\mathbf{H}_{\sigma}\right) \quad|\sigma(t, x)-\sigma(t, y)|^{2} \leq \lambda(t)|x-y|^{2} \kappa_{3}(|x-y|),
	\end{eqnarray*}
	where $\lambda(t)>0$  is locally square integrable and $\kappa_{i} $ is a positive continous function, bounded on $[1, \infty) $ and satisfying
	$$
	\lim _{x \downarrow 0} \frac{\kappa_{i}(x)}{\log x^{-1}}=\rho_{i}<\infty, \quad i=1,2,3.
	$$
	\item [(H2)] There exists a positive constant $L_1$ such that
	\begin{eqnarray*}
		\left(\mathbf{H}_{f,g,\sigma}\right) \quad|f(t, x)|^2+|g(t, x)|^2 +|\sigma(t, x)|^{2} \leq L_1(1+|x|^2).
	\end{eqnarray*}
\end{enumerate}

Comparing our conditions with (H0), one will find that in our conditions the modulus of continuity for $f,g$ in $x$ is different from that for $\sigma$, which is convenient to control. Besides, $f,g, \sigma$ do not depend
on all the path but only on the value at $t$. Therefore, our conditions are more general. For example, consider
$$
f(t, x) :=\lambda(t) \sum_{k \geqslant 1} \frac{\sin (k x)}{k^{2}},
$$
where $\lambda(t)$ is continuous, bounded on $(0, 1]$ and locally square integrable. Then by
Lemma 3.1 and Lemma 4.1 in \cite{Airault2002} and the H\"{o}lder's inequality, we have
$$
\begin{aligned}|f(t, x)-f(t, y)| & \leq \lambda(t) \sum_{k \geqslant 1} \frac{|\sin (k x)-\sin (k y)|}{k^{2}} \\ & \leq 2 \lambda(t) \sum_{k \geqslant 1} \frac{|\sin (k(x-y)) / 2|}{k^{2}} \\ & \leq C \lambda(t)|x-y| \tilde{\kappa}_{1}(|x-y|), \end{aligned}
$$
where
$$
\tilde{\kappa}_{1}(x) :=\left\{\begin{array}{ll}{\log x^{-1},} & {0<x \leq \eta,} \\ {\log \eta^{-1}-1+\frac{\eta}{x},} & {x>\eta.}\end{array}\right.
$$
We take $
\kappa_{1}(x) :=C \tilde{\kappa}_{1}(x)
$, $C$ is a constant. It is easily justified that $\kappa_{1}(x)$
satisfy (H1).
Note that $f(t, x)$ does not satisfy the condition (H0) because $$
\log x^{-1}<(\log x^{-1})^{2},
$$ for $0<x \leq \eta$. Thus our conditions are more general in some sense.

Thus, in order to fill these gaps, assume that $\partial\varphi$ is the subdifferential operator associated to $\varphi$, we will study the averaging principle of Eq. (\ref{orig1}) in sublinear expectation space under non-Lipschitz condition (H1) which is weaker than the classical non-Lipschitz condition (H0) and Lipschitz condition when one discusses variable applications in real world.

This work is organized as follows. Section 2 introduces some notations
and preliminaries which will be useful in what follows. Section 3 is devoted
to established an averaging principle for MSDEs driven by $G$-Brownian motion under non-Lipschitz condition. One example is presented in Section 4 to illustrate our theory.
\section{Preliminaries}
We firstly recall the notion of sublinear expectation space $(\Omega,\mathcal{H},\hat{\mathbb{E}})$, $G$-normal distributed, $G$-Brownian motion and $G$-expectation, see e.g. \cite{Peng2007,Peng2008a,Peng2008b}.
\begin{definition}{\rm
		Let $\Omega$ be a given set and let $\mathcal{H}$ be a vector lattice of real valued functions defined on $\Omega$, namely $c\in \mathcal{H}$ for each constant $c$ and $|X|\in \mathcal{H}$, if $X\in \mathcal{H}$. $\mathcal{H}$ is considered as the space of random variables. A sublinear expectation $\hat{\mathbb{E}}$ on $\mathcal{H}$ is a function $\hat{\mathbb{E}}:\mathcal{H}\rightarrow \mathbb{R}$ satisfying the following properties: for all $X,Y\in \mathcal{H},$ we have
		\begin{description}
			\item[(i)] monotonicity: if $X\geq Y$, then $\hat{\mathbb{E}}[X]\geq\hat{\mathbb{E}}[Y]$;
			\item[(ii)] constant preserving: $\hat{\mathbb{E}}[c]=c$;
			\item[(iii)] sub-additivity: $\hat{\mathbb{E}}[X+Y]\leq \hat{\mathbb{E}}[X]+\hat{\mathbb{E}}[Y]$;
			\item[(iv)] positive homogeneity: $\hat{\mathbb{E}}[\lambda X]=\lambda\hat{\mathbb{E}}[X]$ for each $\lambda\geq0$.
		\end{description}
		
		The triple $(\Omega,\mathcal{H},\hat{\mathbb{E}})$ is called a sublinear expectation space. $X\in \mathcal{H}$ is called a random variable in $(\Omega,\mathcal{H},\hat{\mathbb{E}})$. }
\end{definition}

\begin{definition}{\rm
		($G$-Normal Distribution). A random variable $X$ on a sublinear expectation space $(\Omega,\mathcal{H},\hat{\mathbb{E}})$ is called (centralized) $G$-normal distributed if for any $a,b\geq 0$
		$aX+b\bar{X}\overset{d}{=}\sqrt{a^2+b^2}X,$
		where $\bar{X}$ is an independent copy of $X$. The letter $G$ denotes the function
		$$G(a):=\frac{1}{2}(\overline{\sigma}^2a^{+}-\underline{\sigma}^2a^{-}),a\in \mathbb{R},$$
		with $\underline{\sigma}^2:=-\hat{\mathbb{E}}[-X^2]\leq\hat{\mathbb{E}}[X^2]=:\overline{\sigma}^2.$}
\end{definition}

\begin{definition}{\rm
		A process $(B_t\in \mathcal{H})_{t\geq 0}$ on a sublinear expectation space $(\Omega,\mathcal{H},\hat{\mathbb{E}})$ is called a $G$-Brownian motion if the following properties are satisfied:
		\begin{description}
			\item[(a)] $B_0=0,$
			\item[(b)]For each $t, s \geq 0$, the increment $B_{t+s}-B_{t}\overset{d}{=}\sqrt{s}X$ is independent from\\ $(B_{t_1},B_{t_2}, \cdots,B_{t_n})$ for each $n \in \mathbb{N}, 0 \leq t_1 \leq t_2 \leq \cdots t_n \leq t$, where $X$ is $G$-normal distributed,
			\item[(c)] $\lim_{t\downarrow 0}\mathbb{E}[|B_t|^3]t^{-1}=0.$
	\end{description}}
\end{definition}

Let $\Omega=C_{0}([0,\infty),\mathbb{R}^{d})$, the space of all $\mathbb{R}^{d}$-valued continuous functions on $[0,\infty)$ with $\omega_{0}=0$, be equipped with the distance
$$\rho( {{\omega ^1},{\omega ^2}} ) = \sum_{i = 1}^\infty  {{2^{ - i}}} [ {( {{{\max }_{t \in [ {0,i}]}}}| {\omega _t^1 - \omega _t^2}|)\wedge 1}],$$
and for every $\omega \in \Omega$, let $B_{t}(\omega)=\omega_{t}$ be the canonical process.

For each $T>0$, denote
\begin{eqnarray*}
	Lip\left( {{\Omega _T}} \right)&:=& \{ \psi (B_{{t_1}}, \ldots ,{B_{{t_n}}}):n \ge 1,{t_{1}} \ldots {t_n} \in [ {0,T}], \psi  \in {C_{b.Lip}}( \mathbb{R}^{d \times n})\},
\end{eqnarray*}
where $C_{b,Lip}(\mathbb{R}^{d \times n})$ is the space of all bounded Lipschitz functions defined on $\mathbb{R}^{n}$ and $Lip\left( \Omega  \right) = \mathop  \cup \limits_T Lip\left( {{\Omega _T}} \right).$

Denote by $\mathbb{S}_d$ the collection of all $d\times d$ symmetric matrices. For each given monotonic and sublinear function $G:\mathbb{S}_d\rightarrow \mathbb{R},$ we can construct a $G$-expectation $\hat{\mathbb{E}}$. In the following, we want to construct a sublinear expectation on $(\Omega,Lip(\Omega))$, such that the canonical process $\{B_t\}_{t\geq 0}$ is  a $G$-Brownian motion. For this we first construct a sequence of $d$-dimensional random vectors $(\xi_i)_{i=1}^{\infty}$ on a sublinear expectation space $(\tilde{\Omega},\tilde{\mathcal{H}},\tilde{\mathbb{E}})$ such that $\xi_i$ is $G$-normal distributed and $\xi_{i+1}$ is independent from $(\xi_1,\cdots \xi_i)$ for each $i=1,2,\cdots$.

We now introduce a sublinear expectation $\hat{\mathbb{E}}$ defined on $Lip(\Omega)$ via the following procedure: for each $X\in Lip(\Omega)$ with
$$X=\psi(B_{t_1}-B_{t_0},B_{t_2}-B_{t_1},\cdots,B_{t_n}-B_{t_{n-1}}),$$
for some $ \psi \in C_{b,Lip}(\mathbb{R}^{d\times n})$ and $0=t_0<t_1<\cdots<t_n<\infty$, we set
\begin{eqnarray*}
	&&\hat{\mathbb{E}}[\psi (B_{t_1}-B_{t_0},B_{t_2}-B_{t_1},\cdots,B_{t_n}-B_{t_{n-1}})]\cr
	&&\quad =:\tilde{\mathbb{E}}[\psi(\sqrt{t_1-t_0}\xi_1,\cdots,\sqrt{t_{n}-t_{n-1}}\xi_n)].
\end{eqnarray*}

The related conditional expectation of $$X=\psi(B_{t_1}-B_{t_0},B_{t_2}-B_{t_1},\cdots,B_{t_n}-B_{t_{n-1}}),$$ under $\Omega_{t_j}$ is defined by
\begin{eqnarray*}
	\hat{\mathbb{E}}[X|\Omega_{t_j}]&=&\hat{\mathbb{E}}[\psi(B_{t_1},B_{t_2}-B_{t_1},\cdots,B_{t_n}-B_{t_{n-1}})|\Omega_{t_j}]\cr
	&=:&\psi(B_{t_1},\cdots,B_{t_{j-1}}),
\end{eqnarray*}
where $$\psi(x_1,\cdots,x_j)=\tilde{\mathbb{E}}[\psi(x_1,\cdots,x_j,\sqrt{t_{j+1}-t_{j}}\xi_{j+1},\cdots,\sqrt{t_{n}-t_{n-1}}\xi_n)].$$

It is easy to check that $\hat{\mathbb{E}}[\cdot]$ consistently defines a sublinear expectation on $Lip(\Omega)$ and $(B_t)_{t\geq 0}$ is a $G$-Brownian motion. Since $Lip(\Omega_T)\subset Lip(\Omega)$, $\hat{\mathbb{E}}$ is also a sublinear expectation on  $Lip(\Omega_T)$. We denote by $L^{p}_{G}(\Omega_T)$ the completion of $Lip(\Omega_T)$ under the norm $(\hat{\mathbb{E}}[|\cdot|^{p}])^{\frac{1}{p}}$, for $p\geq 1$.
Therefore, it can
be extended continuously to the completion $L^{p}_{G}(\Omega_T)$.

Let $\mathcal{P}$ be a weakly compact set that represents $\hat{\mathbb{E}}$. For this $\mathcal{P}$, we define the capacity
$$\overline{\mathbb{C}}(A):=\sup_{P\in\mathcal{P}}P(A),A\in\mathcal{B}(\Omega_T).$$
\begin{lemma}\label{lemcapi}
	(see e.g. \cite{Peng2008b}) Let $X\in L_{G}^p$. Then, for each $\alpha>0$, it holds that
	$$\overline{\mathbb{C}}\{|X|>\alpha\}\leq\frac{\hat{\mathbb{E}}[|X|^{p}]}{\alpha^p}.$$
\end{lemma}

Peng (see e.g. \cite{Peng2008b}) also introduced the related stochastic calculus of It\^{o}'s type with respect to $G$-Brownian
motion.
Let $T\in \mathbb{R}_{+}$ be fixed.
\begin{definition}{\rm
		For fixed $p\geq 1, T\in \mathbb{R}_{+}$, define the space $M_{G}^{p,0}([0,T])$ of simple process by
		\begin{eqnarray*}
			M_G^{p,0}([ 0,T ]): = \{\eta _t(\omega): = \sum\limits_{j = 0}^{N - 1} \zeta _{t_j}(\omega) I_{[t_j,t_{j + 1})};\zeta _{t_j}(\omega)  \in L_G^p (\Omega _{t_j}),\\
			\forall N \geq 1,0 = t_0 < t_1 <  \cdots  < t_N = T,j = 0,1, \cdots ,N - 1\}.
		\end{eqnarray*}
		Denote by $M_G^{p}\left( {\left[ {0,T} \right]} \right)$ the completion of $M_G^{p,0}\left( {\left[ {0,T} \right]} \right)$ under the norm
		$M_G^{p,0}\left( {\left[ {0,T} \right]} \right)$  under the norm
		$$\left| \eta  \right|_{M_G^p\left( {\left[ {0,T} \right]} \right)} = \bigg(\int_0^T  \hat{\mathbb{E}}[|\eta_ {s}|^p]  ds  \bigg)^{\frac{1}{p}}.$$}
\end{definition}

\begin{definition} {\rm
		For ${\eta _t}\left( \omega  \right) = \sum\limits_{j = 1}^{N - 1} {{\zeta _{{t_j}}}} \left( \omega  \right){I_{\left[ {{t_j},{t_{j + 1}}} \right)}} \in M_G^{p,0}\left( {\left[ {0,T} \right]} \right)$, define
		$$I\left( \eta  \right) = \int_0^T {{\eta _t}d{B_t}: = \sum\limits_{j = 0}^{N - 1} {{\zeta _j}} } \left( {{B^N_{{t_{j + 1}}}} - {B^N_{{t_j}}}} \right).$$
		The map $I: M_G^{2,0}\left( {\left[ {0,T} \right]} \right) \to L_G^{2}\left( {{\Omega _T}} \right)$ is linear and continuous. Hence, it can be extended continuously to $I: M_G^{2}\left( {\left[ {0,T} \right]} \right) \to L_G^{2}\left( {{\Omega _T}} \right)$. For each $\eta\in M_G^{p}\left( {\left[ {0,T} \right]} \right)$, the stochastic integral is defined by
		$$I( \eta ) = \int_0^T {\eta _t}d{B_t},\eta_t\in M_G^{p}\left( {\left[ {0,T} \right]} \right).$$}
\end{definition}

Unlike the classical theory, the quadratic variation process of $G$-Brownian motion $B$ is not
always a deterministic process and it can be formulated in $ L_G^{2}\left( {{\Omega _t}} \right)$ by
$$\langle B \rangle_t=\lim_{N \rightarrow \infty}\sum_{i=0}^{N-1}(B_{t_{i+1}^N}-B_{t_{i}^N})^2=B_t^2-2\int_{0}^{t}B_sdB_s,$$
where $t_i^N=\frac{iT}{N}$ for each integer $N\geq 1$.
\begin{definition}{\rm
		Define a mapping $M_G^{1,0}\left( {\left[ {0,T} \right]} \right) \to L_G^{1}\left( {{\Omega _T}} \right)$:
		$$Q(\eta)=\int_0^T {\eta _s}\langle B \rangle_s:=\sum_{i=0}^{N-1}(\langle B \rangle_{t_{i+1}^N}-\langle B \rangle_{t_{i}^N})^2,$$
		Then $Q$ can be uniquely extended to $M_G^{1}\left( {\left[ {0,T} \right]} \right) \to L_G^{1}\left( {{\Omega _T}} \right)$. We also denote this mapping by
		$$Q(\eta)=\int_0^T {\eta _s}\langle B \rangle_s.$$}
\end{definition}

In order to prove the main result, we recall the following B-D-G type inequalities, in view of the properties of $G$-Brownian motion and the quadratic
variation process $\langle B \rangle$, see e.g., Gao \cite{Gao2009}, Song \cite{Song2011} and Bai \cite{Bai2014}.

\begin{lemma}
	Let $p\geq 1,\eta \in M_G^{p}\left( {\left[ {0,T} \right]} \right)$ and $0 \leq t \leq T$. Then, we have
	$$\hat{\mathbb{E}}\left[\sup_{0\leq t \leq T}\bigg|\int_{0}^{t}\eta(s)d\langle B\rangle_{s}\bigg|^{p}\right]\leq T^{p-1} \hat{\mathbb{E}}\left
	[\int_{0}^{T}|\eta(s)|^{p}ds\right].$$
\end{lemma}

\begin{lemma}\label{lemBDG}
	Let $p\geq 1,\eta \in M_G^{p}\left( {\left[ {0,T} \right]} \right)$ and $0 \leq t \leq T$. Then, for $C_{p,T}>0$, we have
	$$\hat{\mathbb{E}}\left[\sup_{0\leq t \leq T}\bigg|\int_{0}^{t}\eta(s)d B_s\bigg|^{p}\right]\leq C_{p,T} \hat{\mathbb{E}}\left[\bigg(\int_{0}^{T}|\eta(s)|^{2}d{s}\bigg)^{\frac{p}{2}}\right].$$
\end{lemma}

Note that the letter $C$ below with or without subscripts will denote
		positive constants whose value may change in different occasions. We will
		write the dependence of constants on parameters explicitly if it is essential.

\section{An averaging principle for non-Lipschitz MSDEs driven by $G$-Brownian motion}
\subsection{Notations and Assumptions}

Now, we consider the following assumptions on the function $\varphi$:
\begin{enumerate}
	\item [(H3)] The function $\varphi:\mathbb{R}^d\rightarrow(-\infty,+\infty)$ is l.s.c. such that $${\rm Int} ({\rm Dom}(\varphi))\neq \emptyset, $$ where
	${\rm Dom}(\varphi) \equiv \{x\in \mathbb{R}^d:{\rm Dom}(\varphi)<+\infty\}$
	and supposes that
	$0\in {\rm Int} ({\rm Dom}(\varphi))$ and $\varphi(x)\geq \varphi(0)=0$, for all $x\in \mathbb{R}^d$.
\end{enumerate}

We recall that the subdifferential of the function $\varphi$ is defined by
$$\partial\varphi(x)=\{y\in \mathbb{R}^d:\langle y,z-x\rangle+\varphi(x)\leq \varphi(z), \vee z\in \mathbb{R}^d\},$$
and by $(x,x^{*})\in Gr(\partial\varphi)$, which means that $x\in {\rm Dom(\partial\varphi)}$ and $x^{*}\in {\rm Dom(\partial\varphi)}$ where
$${\rm Dom(\partial\varphi)}:=\{x\in \mathbb{R}^d:\partial\varphi(x)\neq \emptyset\}.$$
\begin{remark}{\rm
		Under assumption (H3), the subdifferential
		operator $\varphi$ becomes a maximal monotone operator, i.e.
		maximal in the class of operators which satisfy the condition
		$$\langle y^{*}-z^{*},y-z\rangle \geq 0, (y,y^{*}),(z,z^{*})\in Gr(\partial\varphi).$$
		Conversely (only in the case $d=1$), we recall that, if $A$ is a given maximal monotone operator on $\mathbb{R}$, then there
		exists a proper l.s.c. function, such that $A=\partial\varphi$.}
\end{remark}

Now we give the previous definition of the solution to Eq. (\ref{orig1}).
\begin{definition}\label{defsolu}
	{\rm The solution of Eq. (\ref{orig1}) is a pair of processes $(X,K)$ satisfying that
		\begin{enumerate}
			\item [(i)] $X(0)=\xi$ and $X(t) \in \overline{{\rm Dom}(\varphi)} $ for every $t \in[s,T],q.s.$;
			\item [(ii)] for any $s<T<+\infty$, $K$ is continuous and of finite variation on $[s,T]$ with $K(0)=0$, and for all $t,t'\in [s,T],$ $$\lim_{t\rightarrow t'}{\hat{\mathbb{E}}}[|K(t)-K(t')|]^2=0;$$
			\item [(iii)]
			\begin{eqnarray*}
				X(t)+K(t)&=&X(s)+\int_{s}^{t}f(r,X(r))dr+ \int_{s}^{t}\sigma(r,X(r))dB_r\cr
				&&+\int_{s}^{t}g(r,X(r))d\langle B\rangle_r, s\leq t \leq T, q.s.
			\end{eqnarray*}
			\item [(iv)]
			$\int_{s}^{t}\langle u-X(r),dK(r)\rangle-\int_{s}^{t} \varphi(X(r))dr \leq (t-s)\varphi(u), u\in \mathbb{R}^d, s\leq t \leq T, q.s.$
	\end{enumerate}}
\end{definition}
\begin{lemma}\label{lem2.4} We shall write $dK(t)\in \varphi(X(t))dt$, q.s., instead of inequality (iv) (see also the below result). Now, taking the processes $X,X',K,K'$ such
	that $dK(t)\in \varphi(X(t))dt$ and $dK'(t)\in \varphi(X'(t))dt$, we
	have
	\begin{eqnarray*}
		\int_{0}^{t}\langle X(t)-X'(t), dK(t)-dK'(t)\rangle \geq 0.
	\end{eqnarray*}
\end{lemma}
\begin{lemma}\label{lem2.5}
	Let $\Phi \in C^2([0,T]\times\mathbb{R})$ be a real function with $\partial_t \Phi,\partial_x\Phi,\partial_{xx} \Phi \in C_{b,Lip(\mathbb{R})}$. Let $f,\sigma$ and $g$ be
	bounded processes in $M^2_G([0,T])$ and $K\in M^2_G([0,T])$ is
	continuous, which satisfies for each $s \leq t \leq T$,
	$$\lim_{t\rightarrow s}{\hat{\mathbb{E}}}[|K(t)-K(s)|^2]=0;$$
	Then, we have
	\begin{eqnarray*}
		\Phi(t, X(t))-\Phi(s, X(s))&=&\int_{s}^{t}\bigg[\frac{\partial \Phi}{\partial t}(X(r))+\frac{\partial \Phi}{\partial x}(X(r))f(r)\bigg]dr\cr
		&&+\int_{s}^{t}\frac{\partial \Phi}{\partial t}(X(r))g(r)d\langle B \rangle_r +\int_{s}^{t}\frac{\partial \Phi}{\partial x}(X(r))\sigma(r)dB_r\cr
		&&-\int_{s}^{t}\frac{\partial \Phi}{\partial x}(X(r))dK(r)+\frac{1}{2}\int_{s}^{t}\frac{\partial^{2} \Phi}{\partial x^2}(X(r))\sigma^2(r)d\langle B \rangle_r.
	\end{eqnarray*}
\end{lemma}

The following generalization of the Gronwall-Belmman type inequality comes from Bihari \cite{Bihari1975}.
\begin{lemma}\label{lembiha}
	Let $h(s)$ be a strictly positive function on $\mathbb{R}_+$ satisfying for some $\delta>0$
	$$
	h(t) \leq h(0)+\delta \int_{0}^{t} \rho_{\eta}(h(s)) \mathrm{d} s, \quad t \geqslant 0,
	$$
	where
	$$
	\rho_{\eta}(x) :=\left\{\begin{array}{ll}{x \log x^{-1},} & {x \leq \eta,} \\ {\eta \log \eta^{-1}+\left(\log \eta^{-1}-1\right)(x-\eta),} & {x>\eta,}\end{array}\right.
	$$
	for $
	0<\eta<1 / e$.
	Then for any  $T>0$, there is a constant $ C:=C(T, \delta, \eta)$ such that
	$$
	h(t) \leq C((h(0))^{\exp \{-\delta T\}}+h(0)), \quad t \in[0, T].
	$$
\end{lemma}
\para{Proof:} This can be found from Zhang's work \cite[Lemma 2.1]{Zhang2005}.
\qed

\subsection{The Unique Solution of Eq. (\ref{orig1}) }
\begin{lemma}
	Assume that (H1)-(H3) hold. Then, Eq. (\ref{orig1}) has a unique solution.
\end{lemma}
\para{Proof:}
The existence and uniqueness of the solution for Eq. (\ref{orig1})  can be shown by means of the penalisation method. We can use the Yosida approximation of the operator $\partial_{\varphi}$ to complete the proof. That is, for $\epsilon\in(0,1]$, let $\nabla \varphi_{\epsilon}$ be the gradient of $\partial_{\varphi}$, where $\partial_{\varphi}$ is the Moreau-Yosida regularisation of $\varphi$, i.e.
$$
\varphi_{\epsilon}(x) :=\inf \left\{\frac{1}{2 \epsilon}|v-x|^{2}+\varphi(v) : v \in \mathbb{R}^{d}\right\}, \quad \epsilon>0,
$$
where $\varphi$, is a $C^1$ convex function. So, on the basis of Theorem 3.1 in \cite{Ren2017} and Theorem 1.2 in \cite{Qiao2014}, there exists the unique solution $X(t), t \in [0,T]$ to Eq. (\ref{orig1}). \qed

\begin{lemma}\label{lemb}
	Suppose that original Eq. {\rm (\ref{orig1})} satisfies the conditions {\rm (H1)-(H3)}. Then, for $p\geq 1$, we have
	$$	\hat{\mathbb{E}}\bigg[\sup_{t\in[0,T]}|X(t)|^{2p}\bigg]\leq C.$$
\end{lemma}
\para{Proof:} By It\^{o}'s formula, we have
\begin{eqnarray*}
	|X(t)|^{2p} &=&|X(0)|^{2p}+2p\int_0^t |X(s)|^{2p-2}\langle X(s), f(s, X(s))\rangle  ds\cr
	&&
	+ 2p  \int_0^t|X(s)|^{2p-2} \langle X(s),g(s, X(s))\rangle  d\langle B\rangle_s\cr
	&&
	+ 2p \int_0^t|X(s)|^{2p-2} \langle X(s),\sigma(s, X(s))\rangle  dB_s\cr
	&&+2p(p-1)  \int_0^t|X(s)|^{2p-4}|X^T(s)\sigma(s, X(s))|^2d\langle B\rangle_s\cr
	&&+p  \int_0^t|X(s)|^{2p-2} |\sigma(s, X(s))|^2d\langle B\rangle_s
	\cr
	&&-2p \int_0^t |X(s)|^{2p-2} \langle X(s),d K(s)\rangle.
\end{eqnarray*}
Firstly, by Lemma \ref{lem2.4}, we have  $$\int_0^t |X(s)|^{2p-2} \langle X(s),d K(s)\rangle\geq0.$$
Then, we have
\begin{eqnarray*}
	\hat{\mathbb{E}}\bigg[\sup_{t\in[0,T]}|X(t)|^{2p}\bigg]&\leq& \hat{\mathbb{E}}[|X(0)|^{2p}]+C \hat{\mathbb{E}}\bigg[\sup_{t\in[0,T]}\bigg|\int_0^t |X(s)|^{2p-2}\langle X(s), f(s, X(s))\rangle  ds\bigg|\bigg]\cr
	&&
	+ C \hat{\mathbb{E}}\bigg[\sup_{t\in[0,T]}\bigg|\int_0^t|X(s)|^{2p-2} \langle X(s),g(s, X(s))\rangle  d\langle B\rangle_s\bigg|\bigg]\cr
	&&
	+ C \hat{\mathbb{E}}\bigg[\sup_{t\in[0,T]}\bigg|\int_0^t|X(s)|^{2p-2} \langle X(s),\sigma(s, X(s))\rangle  dB_s\bigg|\bigg]\cr
	&&+C  \hat{\mathbb{E}}\bigg[\sup_{t\in[0,T]}\bigg|\int_0^t|X(s)|^{2p-4}|X^T(s)\sigma(s, X(s))|^2d\langle B\rangle_s\bigg|\bigg]\cr
	&&+C \hat{\mathbb{E}}\bigg[\sup_{t\in[0,T]}\bigg|\int_0^t|X(s)|^{2p-2} |\sigma(s, X(s))|^2d\langle B\rangle_s\bigg|\bigg]\cr
	&\leq& \hat{\mathbb{E}}[|X(0)|^{2p}]+C \int_0^T \hat{\mathbb{E}}[|X(s)|^{2p-1}|f(s, X(s))|]  ds\cr
	&&
	+C \int_0^T \hat{\mathbb{E}}[|X(s)|^{2p-1}|g(s, X(s))|]  ds\cr
	&&+C \int_0^T \hat{\mathbb{E}}[|X(s)|^{2p-2}|\sigma(s, X(s))|^2]  ds\cr
	&&
	+ C \hat{\mathbb{E}}\bigg[\sup_{t\in[0,T]}\bigg|\int_0^t|X(s)|^{2p-2} \langle X(s),\sigma(s, X(s))\rangle dB_s\bigg|\bigg].
\end{eqnarray*}
Next, by Lemma \ref{lemBDG}, we have
\begin{eqnarray*}
	&&\hat{\mathbb{E}}\bigg[\sup_{t\in[0,T]}\bigg|\int_0^t|X(s)|^{2p-2} \langle X(s),\sigma(s, X(s))\rangle  dB_s\bigg|\bigg]\cr
	&&\leq C \hat{\mathbb{E}}\bigg[\int_0^T|X(s)|^{4p-2}|\sigma(s, X(s))|^2d\langle B\rangle_s\bigg]^{\frac{1}{2}}\cr
	&&\leq \frac{1}{2p} \hat{\mathbb{E}}\bigg[\sup_{t\in[0,T]}|X(t)|^{2p}\bigg]+C\int_0^T\hat{\mathbb{E}}[|X(s)|^{2p-2}|
	\sigma(s, X(s))|^{2}]ds.
\end{eqnarray*}
Finally, we have
\begin{eqnarray*}
	\hat{\mathbb{E}}\bigg[\sup_{t\in[0,T]}|X(t)|^{2p}\bigg]&\leq&  C(1+\hat{\mathbb{E}}[|X(0)|^{2p}])+C \int_0^T \hat{\mathbb{E}}[\sup_{s\in[0,t]}|X(s)|^{2p}] dt \cr &&
	+ \frac{1}{2p} \hat{\mathbb{E}}\bigg[\sup_{t\in[0,T]}|X(t)|^{2p}\bigg].
\end{eqnarray*}
Therefore, by Gronwall's inequality, we have
\begin{eqnarray*}
	\hat{\mathbb{E}}\bigg[\sup_{t\in[0,T]}|X(t)|^{2p}\bigg]\leq C.
\end{eqnarray*}
This completed the proof. \qed

\subsection{The Averaging Principle}
Consider the following MSDEs driven by $G$-Brownian motion on $\mathbb{R}^d$:
\begin{eqnarray}\label{orginte1}
X^\varepsilon(t)&=&X(0)+\int_{0}^{t}f(s/\varepsilon,X^\varepsilon(s))ds+ \int_{0}^{t}\sigma(s/\varepsilon,X^\varepsilon(s))dB_s
\cr&&
+\int_{0}^{t}g(s/\varepsilon,X^\varepsilon(s))d\langle B\rangle_s-K(t), t \in [0,T],
\end{eqnarray}
where $ \varepsilon \in
(0,\varepsilon_0]$ is a positive small parameter with $ \varepsilon_0 $ a
fixed number.
This paper will study convergence $X^{\varepsilon}(t)\rightarrow \bar{X}(t),\varepsilon\rightarrow 0$, where $\bar{X}(t)$ is the solution of the averaged equation
\begin{eqnarray}\label{orginte2}
\bar{X}(t)&=&X(0)+\int_{0}^{t}\bar f(\bar{X}(s))ds+ \int_{0}^{t}\bar{\sigma}(\bar{X}(s))dB_s
\cr&&
+\int_{0}^{t}\bar g(\bar{X}(s))d\langle B\rangle_s-\bar K(t),t \in [0,T].
\end{eqnarray}
Under the similar conditions such as
$ X(t)$ in Eq. (\ref{orig1}), this equation will have a unique solution $ \bar{X}(t) $.

To proceed,  with the natural time scaling $t\rightarrow t\varepsilon$, the equivalent form of Eq.(\ref{orginte1})  on $t\in[0,T]$ can be rewritten as the following form on $t\in [0, T/\varepsilon]$.
\begin{eqnarray}\label{orginteg3}
X_\varepsilon(t)&=&X(0)+\varepsilon \int_0^t f(s, X_\varepsilon(s)) ds
+ \sqrt{\varepsilon}  \int_0^t \sigma(s, X_\varepsilon(s))
dB_s\cr&&
+{\varepsilon}\int_{0}^{t}g(s, X_\varepsilon(s))d\langle B\rangle_s- \varepsilon K(t),
\end{eqnarray}
where the coefficients have the
same conditions as in Eq. (\ref{orginte1}).

Similarly, with the natural time scaling, the equivalent form of averaged Eq.(\ref{orginte2})  on $t\in[0,T]$ can be rewritten as the following form on $t\in [0, T/\varepsilon]$:
\begin{eqnarray}\label{ave}
Z_\varepsilon(t)&=&X(0)+\varepsilon \int_0^t \bar{f} (Z_\varepsilon(s)) ds
+ \sqrt{\varepsilon}  \int_0^t \bar{\sigma}  (Z_\varepsilon(s))
dB_s
\cr&&
+{\varepsilon}\int_{0}^{t}\bar{g}(Z_\varepsilon(s))d\langle B\rangle_s-\varepsilon  \bar{K}(t).
\end{eqnarray}

Thus, the study of convergence $X^{\varepsilon}(t)\rightarrow \bar{X}(t),\varepsilon\rightarrow 0$ on finite time intervals is equivalent to the study of convergence $X_\varepsilon(t)\rightarrow Z_\varepsilon(t) ,\varepsilon\rightarrow 0$ in time intervals of order $\varepsilon^{-1}$. So we will claim the main theorem to show relationship between solution processes $ X_\varepsilon(t) $ to
the orignal Eq. (\ref{orginteg3})  and $ Z_\varepsilon(t) $ to the averaged Eq. (\ref{ave}) in time intervals of order $\varepsilon^{-1}$. It shows that the solution of averaged Eq. (\ref{ave}) converges to that
of the original Eq. (\ref{orginteg3}) in the sense of $p$-th moments and convergence in capacity.

In order to establish the averaging principle, we need assume further that  the functions $\bar{f}: \mathbb{R}^d \rightarrow
\mathbb{R}^d,\bar{g}$, $\bar{\sigma}: \mathbb{R}^d \rightarrow \mathbb{R}^{d \times n}$ are continuous. Presuming they meet the following additional inequalities:
\begin{eqnarray*}
	&&\left(\mathbf{C}_{f}\right) \quad\frac{1}{T_1}\int_0^{T_1} | f(s,x)-\bar{f}(x)
	| ds
	\leq \varphi_1(T_1)(1+ |x|),\cr
	&&\left(\mathbf{C}_{g}\right) \quad\frac{1}{T_1}\int_0^{T_1} | g (s,x)-\bar{g}(x)| ds
	\leq \varphi_2(T_1)(1+ |x|),\cr
	&&\left(\mathbf{C}_{\sigma}\right) \quad\dfrac{1}{T_1} \int_0^{T_1} |\sigma(s,x)-\bar{\sigma }(x)|^2ds
	\leq \varphi_3(T_1)(1+|x|^2),
\end{eqnarray*}
where $ T_1 \in [0,T],  \varphi_i(T_1)$ are positive bounded functions
with $$\lim_{T_1 \rightarrow \infty} \varphi_i(T_1)= 0, i = 1,2,3.$$

\begin{theorem}\label{thm1}
	Suppose that original Eq. (\ref{orginteg3}) and averaged Eq. (\ref{ave})  both satisfy the assumptions {\rm (H1)-(H3)} and {\rm $(\mathbf{C}_{f}), (\mathbf{C}_{g}), (\mathbf{C}_{\sigma})$}. For a given arbitrarily small number $ \delta_1 > 0,$ there exist $ L > 0 $, $\alpha \in (0,1) $, $ \varepsilon_1 \in (0,\varepsilon_0] $, such that for any $ \varepsilon \in (0,\varepsilon_1] $ , each $ t \in
	[0,L\varepsilon^{\frac{1}{2}-\alpha}] $, we have
	\begin{eqnarray*}
		\hat{\mathbb{E}}\bigg[\sup_{t\in[0,L\varepsilon^{\frac{1}{2}-\alpha}]}|X_\varepsilon(t)- Z_\varepsilon(t)|^{2p}\bigg]\leq \delta_1.
	\end{eqnarray*}
\end{theorem}
We also have the following result on uniform convergence in capacity.
\begin{corollary}\label{thm2}
	Suppose that all assumptions {\rm (H1)-(H3)} and {\rm$(\mathbf{C}_{f}), (\mathbf{C}_{g}), (\mathbf{C}_{\sigma})$} are satisfied. Then for any number $ \delta_2 > 0 $,  each $ t \in
	[0,L\varepsilon^{\frac{1}{2}-\alpha}] $, we have
	\begin{eqnarray*}
		\lim_{\varepsilon \rightarrow 0}\overline{\mathbb{C}}\bigg[\sup_{t\in[0,L\varepsilon^{\frac{1}{2}-\alpha}]}|X_\varepsilon(t)- Z_\varepsilon(t)|>\delta_2\bigg] = 0,
	\end{eqnarray*}
	where $ L $ and $ \alpha $ are the same to Theorem \ref{thm1}.
\end{corollary}
\para{The Proof of Theorem \ref{thm1}:}
By It\^{o}'s formula, let $\Delta(s)=X_\varepsilon(s)-Z_\varepsilon(s)$, we have
\begin{eqnarray*}
	|X_\varepsilon(t)-Z_\varepsilon(t)|^{2p} &=& 2p\varepsilon \int_0^t |\Delta(s)|^{2p-2}\langle \Delta(s), f(s, X_\varepsilon(s))-\bar{f} (Z_\varepsilon(s))\rangle  ds\cr
	&&
	+ 2p{\varepsilon}  \int_0^t|\Delta(s)|^{2p-2} \langle \Delta(s),g(s, X_\varepsilon(s))-\bar{g} (Z_\varepsilon(s))\rangle  d\langle B\rangle_s\cr
	&&
	+ 2p\sqrt{\varepsilon}  \int_0^t|\Delta(s)|^{2p-2} \langle \Delta(s),\sigma(s, X_\varepsilon(s))-\bar{ \sigma} (Z_\varepsilon(s))\rangle  dB_s\cr
	&&+2p(p-1){\varepsilon}  \int_0^t|\Delta(s)|^{2p-4}|\Delta^{T}(s)(\sigma(s, X_\varepsilon(s))-\bar{ \sigma} (Z_\varepsilon(s)))|^2d\langle B\rangle_s\cr
	&&+p{\varepsilon}  \int_0^t|\Delta(s)|^{2p-2} |\sigma(s, X_\varepsilon(s))-\bar{ \sigma} (Z_\varepsilon(s))|^2d\langle B\rangle_s
	\cr
	&&-2p\varepsilon \int_0^t |\Delta(s)|^{2p-2} \langle X_\varepsilon(s)-Z_\varepsilon(s),
	d K(s)- d\bar{K}(s)\rangle\cr
	&=:&\sum_{i=1}^{6}I_i.
\end{eqnarray*}

Firsly, for $I_1$, by the elementary inequality, we can obtain
\begin{eqnarray*}
	\hat{\mathbb{E}}[\sup_{0 \leq t \leq u}I_1] &\leq&  2p\varepsilon \hat{\mathbb{E}}\bigg[\sup_{0 \leq t \leq u} \bigg|\int_0^t |\Delta(s)|^{2p-2}\langle \Delta(s), f(s, X_\varepsilon(s))-\bar{f} (Z_\varepsilon(s))\rangle  ds\bigg|\bigg]\cr
	&\leq&2p\varepsilon \hat{\mathbb{E}}\bigg[\sup_{0 \leq t \leq u}\bigg|\int_0^t |\Delta(s)|^{2p-2}\langle \Delta(s), f(s, X_\varepsilon(s))-f(s, Z_\varepsilon(s))\rangle  ds\bigg|\bigg] \cr
	&&+2p\varepsilon \hat{\mathbb{E}}\bigg[\sup_{0 \leq t \leq u}\bigg|\int_0^t |\Delta(s)|^{2p-2}\langle \Delta(s),f(s, Z_\varepsilon(s))-\bar{f} (Z_\varepsilon(s))\rangle  ds\bigg|\bigg] \cr
	&
	=:&I_{11}+I_{12}.
\end{eqnarray*}
By (H1), for $I_{11}$, one has
\begin{eqnarray}\label{i11}
I_{11}&=&2p\varepsilon \hat{\mathbb{E}}\bigg[\sup_{0 \leq t \leq u}\bigg|\int_0^t |\Delta(s)|^{2p-2}\langle \Delta(s), f(s, X_\varepsilon(s))-f(s, Z_\varepsilon(s))\rangle  ds\bigg| \bigg]\cr
&&
\leq 2p\varepsilon \int_0^u \lambda(s) \hat{\mathbb{E}} \bigg[|X_\varepsilon(s)-Z_\varepsilon(s)|^{2p}\kappa_1(|X_\varepsilon(s)-Z_\varepsilon(s)|)\bigg]ds.
\end{eqnarray}
Next, for $ I_{12} $, by Young's inequality, and {\rm (C1)}, we get
\begin{eqnarray*}
	I_{12} &\leq&2p\varepsilon \hat{\mathbb{E}}\bigg[\sup_{0 \leq t \leq u}\bigg|\int_0^t |\Delta(s)|^{2p-2}\langle \Delta(s), f(s, Z_\varepsilon(s))-\bar{f}(Z_\varepsilon(s))\rangle ds\bigg|\bigg]\cr
	&\leq&\varepsilon \hat{\mathbb{E}}\bigg[\sup_{0 \leq t \leq u}\bigg|\int_0^t(C_{p}(2p-1)+1) |\Delta(s)|^{2p}|f(s, Z_\varepsilon(s))-\bar{f}(Z_\varepsilon(s))ds\bigg|\bigg]\cr
	&\leq &C_{12}(C_{p}(2p-1)+1)\varepsilon \hat{\mathbb{E}}\bigg[\sup_{0 \leq t \leq u} t\dfrac{1}{t}\int_0^t|f(s, Z_\varepsilon(s))-\bar{f}(Z_\varepsilon(s))| ds\bigg]\cr
	&\leq &C_{12} (C_{p}(2p-1)+1)\varepsilon u [\sup_{0 \leq t \leq u} \varphi_1(t) ]  \bigg(1+\hat{\mathbb{E}}\bigg[\sup_{0 \leq s \leq u} |Z_\varepsilon (s)|^2 \bigg]\bigg).
\end{eqnarray*}
Similar to the proof of Lemma \ref{lemb}, then, we have $$ \hat{\mathbb{E}}[\sup_{0 \leq s \leq u} |Z_\varepsilon (s)|^2]<  \infty. $$ This property combines with the boundedness of $ \varphi_1(t) $. We can further estimate that there exists a constant $ C_{12}$ such that
\begin{eqnarray}\label{i12}
I_{12} \leq  \varepsilon u C_{12}.
\end{eqnarray}
So putting (\ref{i11}) and (\ref{i12}) together, we have
\begin{eqnarray}\label{i1}
I_{1} \leq  2p\varepsilon \int_0^u \lambda(s) \hat{\mathbb{E}} \bigg[|X_\varepsilon(s)-Z_\varepsilon(s)|^{2p}\kappa_1(|X_\varepsilon(s)-Z_\varepsilon(s)|)\bigg]ds
+ \varepsilon u C_{12}.
\end{eqnarray}

Then for $ I_2 $, we get
\begin{eqnarray}\label{i2-}
\hat{\mathbb{E}}[\sup_{0 \leq t \leq u}I_2] &\leq&2p{\varepsilon}  \hat{\mathbb{E}}\bigg[\sup_{0 \leq t \leq u}\bigg|\int_0^t|\Delta(s)|^{2p-2} \langle \Delta(s),g(s, X_\varepsilon(s))-\bar{g} (Z_\varepsilon(s))\rangle  d\langle B\rangle_s\bigg|\bigg]\cr
&\leq&2p{\varepsilon}  \hat{\mathbb{E}}\bigg[\int_0^u|\Delta(s)|^{2p-2} \langle \Delta(s),g(s, X_\varepsilon(s))-\bar{g} (Z_\varepsilon(s))\rangle  ds\bigg]\cr
&\leq&2p\varepsilon \hat{\mathbb{E}}\bigg[ \int_0^u|\Delta(s)|^{2p-2}\langle \Delta(s),g(s, X_\varepsilon(s))-g(s, Z_\varepsilon(s))\rangle ds\bigg]\cr
&&+2p\varepsilon \hat{\mathbb{E}}\bigg[ \int_0^u|\Delta(s)|^{2p-2}\langle \Delta(s),g(s, Z_\varepsilon(s))-\bar{g}  (Z_\varepsilon(s))\rangle                                                                   ds\bigg]\cr
&\leq&2pL\varepsilon \int_0^u \lambda(s) \hat{\mathbb{E}} \bigg[|X_\varepsilon(s)-Z_\varepsilon(s)|^{2p}\kappa_2(|X_\varepsilon(s)-Z_\varepsilon(s)|)\bigg]ds\cr
&&+{\varepsilon}  \hat{\mathbb{E}}\bigg[\int_0^u (C_{p}(2p-1)+1) |\Delta(s)|^{2p})|g(s, X_\varepsilon(s))-\bar{g} (Z_\varepsilon(s))|) ds\bigg]\cr
&\leq&2p\varepsilon \int_0^u \lambda(s) \hat{\mathbb{E}} \bigg[|X_\varepsilon(s)-Z_\varepsilon(s)|^{2p}g_2(|X_\varepsilon(s)-Z_\varepsilon(s)|)\bigg]ds\cr
&&+C_{22}(C_{p}(2p-1)+1)\varepsilon \hat{\mathbb{E}}\bigg[u\dfrac{1}{u}\int_0^u|g(s, Z_\varepsilon(s))-\bar{g}(Z_\varepsilon(s))|ds\bigg]\cr
&\leq &2p\varepsilon \int_0^u \lambda(s) \hat{\mathbb{E}} \bigg[|X_\varepsilon(s)-Z_\varepsilon(s)|^{2p}\kappa_2(|X_\varepsilon(s)-Z_\varepsilon(s)|)\bigg]ds\cr
&&+ C_{22}(C_{p}(2p-1)+1)\varepsilon u  [\sup_{0 \leq t \leq u} \varphi_3(t) ]  \bigg(1+\hat{\mathbb{E}} \bigg[\sup_{0 \leq s \leq u} |Z_\varepsilon (s)|^2\bigg]\bigg)\cr
&\leq& 2p\varepsilon \int_0^u \lambda(s) \hat{\mathbb{E}} \bigg[|X_\varepsilon(s)-Z_\varepsilon(s)|^{2p}\kappa_2(|X_\varepsilon(s)-Z_\varepsilon(s)|)\bigg]ds\cr
&&+ \varepsilon uC_{22}.
\end{eqnarray}

To proceed, using the Young's inequality, we have
\begin{eqnarray*}
	I_{3} &\leq&
	2p\sqrt{\varepsilon}\hat{\mathbb{E}} \left[\sup_{0 \leq t \leq u}\bigg|\int_0^t|\Delta(s)|^{2p-2} \langle \Delta(s),\sigma(s, X_\varepsilon(s))-\bar{ \sigma} (Z_\varepsilon(s))\rangle  dB_s\bigg|\right]\cr
	&\leq & 2p \sqrt{\varepsilon}  \hat{\mathbb{E}} \left[\int_0^u |\Delta(s)|^{4p-4}|\langle X_\varepsilon(s)-Z_\varepsilon(s),\sigma(s, X_\varepsilon(s))-\bar{ \sigma} (Z_\varepsilon(s))\rangle |^2ds\right]^{\frac{1}{2}}\cr
	& \leq & 2p \sqrt{\varepsilon} \hat{\mathbb{E}}\left[\sup_{0 \leq s \leq u}|\Delta(s)|^{2p}\int_0^u|\Delta(s)|^{2p-2}|\sigma(s, X_\varepsilon(s))-\bar{ \sigma} (Z_\varepsilon(s))|^2 ds\right]^{\frac{1}{2}}\cr
	& \leq & \frac{1}{2}\hat{\mathbb{E}}\bigg[\sup_{0 \leq t \leq u}|X_\varepsilon(t)-Z_\varepsilon(t)|^{2p}\bigg]\cr
	&&
	+C_p\sqrt{\varepsilon} \hat{\mathbb{E}}\left[\int_0^u|\Delta(s)|^{2p-2}|\sigma(s, X_\varepsilon(s))-\bar{ \sigma} (Z_\varepsilon(s))|^{2} ds\right]\cr
	&=:& I_{31}+I_{32}.
\end{eqnarray*}
According to (\ref{i2-}), we get
\begin{eqnarray*}
	I_{32}& \leq&  C_{31}\sqrt{\varepsilon} \int_0^u \lambda(s) \hat{\mathbb{E}} \bigg[|X_\varepsilon(s)-Z_\varepsilon(s)|^{2p}\kappa_3(|X_\varepsilon(s)-Z_\varepsilon(s)|)\bigg]ds\cr
	&&+C_{32}(C_{p}(2p-1)+1)\sqrt{\varepsilon} \hat{\mathbb{E}}\bigg[u\dfrac{1}{u}\int_0^u|\sigma(s, Z_\varepsilon(s))-\bar{\sigma}(Z_\varepsilon(s))|^2ds\bigg].
\end{eqnarray*}
Thus, we have
\begin{eqnarray}\label{i3}
\hat{\mathbb{E}}[\sup_{0 \leq t \leq u}I_3] &\leq& \frac{1}{2}\hat{\mathbb{E}} \bigg[\sup_{0 \leq s \leq u}|X_\varepsilon(s)-Z_\varepsilon(s)|^{2p}\bigg]+C_{32}u\sqrt{\varepsilon}\cr
&&+C_{31}\sqrt{\varepsilon} \int_0^u \lambda(s) \hat{\mathbb{E}} \bigg[|X_\varepsilon(s)-Z_\varepsilon(s)|^{2p}\kappa_3(|X_\varepsilon(s)-Z_\varepsilon(s)|)\bigg]ds.
\end{eqnarray}

For the term $ I_4 $, we have
\begin{eqnarray*}
	\hat{\mathbb{E}}[\sup_{0 \leq t \leq u}I_4]&\leq&2p(p-1){\varepsilon}	\hat{\mathbb{E}}\bigg[\sup_{0 \leq t \leq u}\bigg|\int_0^t|\Delta(s)|^{2p-4}|\Delta^{T}(s)(\sigma(s, X_\varepsilon(s))-\bar{ \sigma} (Z_\varepsilon(s)))|^2d\langle B\rangle_s\bigg|\bigg]\cr
	&\leq&2p(p-1){\varepsilon}	\hat{\mathbb{E}}\bigg[\int_0^u|\Delta(s)|^{2p-2}|(\sigma(s, X_\varepsilon(s))-\bar{ \sigma} (Z_\varepsilon(s)))|^2ds\bigg].
\end{eqnarray*}
Similar to the proof of $I_{32}$, we have
\begin{eqnarray*}
	\hat{\mathbb{E}}[\sup_{0 \leq t \leq u}I_4]&\leq&C_{41}{\varepsilon}  \int_0^u \lambda(s) \hat{\mathbb{E}} \bigg[|X_\varepsilon(s)-Z_\varepsilon(s)|^{2p}\kappa_3(|X_\varepsilon(s)-Z_\varepsilon(s)|)\bigg]ds\cr
	&&+C_{42}u{\varepsilon}.
\end{eqnarray*}

Next, for the term $ I_5 $, it is easy to obtain
\begin{eqnarray*}
	\hat{\mathbb{E}}[\sup_{0 \leq t \leq u}I_5]&\leq& p{\varepsilon}  \hat{\mathbb{E}}\bigg[\sup_{0 \leq t \leq u}\bigg|\int_0^t|\Delta(s)|^{2p-2} |\sigma(s, X_\varepsilon(s))-\bar{ \sigma} (Z_\varepsilon(s))|^2d\langle B\rangle_s\bigg|\bigg]
	\cr
	&\leq& p{\varepsilon}  \hat{\mathbb{E}}\bigg[\int_0^u|\Delta(s)|^{2p-2} |\sigma(s, X_\varepsilon(s))-\bar{ \sigma} (Z_\varepsilon(s))|^2ds\bigg]\cr
	&\leq&  C_{51}{\varepsilon} \int_0^u \lambda(s) \hat{\mathbb{E}} \bigg[|X_\varepsilon(s)-Z_\varepsilon(s)|^{2p}\kappa_3(|X_\varepsilon(s)-Z_\varepsilon(s)|)\bigg]ds\cr
	&&+C_{52}u{\varepsilon}.
\end{eqnarray*}

Then, according to (\ref{defsolu}) and (\ref{lem2.4}), it is immediate to conclude that
\begin{eqnarray}\label{k}
\int_0^t |\Delta(s)|^{2p-2} \langle X_\varepsilon(s)-Z_\varepsilon(s),
d K(s)- d\bar{K}(s)\rangle \geqslant 0.
\end{eqnarray}

Finally, we have
\begin{eqnarray*}
	&&	\hat{\mathbb{E}}\bigg[\sup_{0 \leq t \leq u}|X_\varepsilon(t)- Z_\varepsilon(t)|^{2p}\bigg]\cr
	&& \ \ \ \leq C_{11}\varepsilon \int_0^u \lambda(s) \hat{\mathbb{E}} \bigg[|X_\varepsilon(s)-Z_\varepsilon(s)|^{2p}g_1(|X_\varepsilon(s)-Z_\varepsilon(s)|)\bigg]ds\cr
	&& \ \ \ +\varepsilon C_{21}\int_0^u \lambda(s) \hat{\mathbb{E}} \bigg[|X_\varepsilon(s)-Z_\varepsilon(s)|^{2p}g_2(|X_\varepsilon(s)-Z_\varepsilon(s)|)\bigg]d\cr
	&& \ \ \ +\frac{1}{2} \hat{\mathbb{E}} \bigg[\sup_{0 \leq s \leq u}|X_\varepsilon(s)-Z_\varepsilon(s)|^{2p}\bigg]+ \varepsilon u C_{12}+\varepsilon uC_{22}+C_{32}u\sqrt{\varepsilon}\cr
	&& \ \ \ +[C_{31}+\sqrt{\varepsilon}(C_{41}+C_{51})]\sqrt{\varepsilon}]\cr
	&& \ \ \  \times \int_0^u \lambda(s) \hat{\mathbb{E}} \bigg[|X_\varepsilon(s)-Z_\varepsilon(s)|^{2p}g_3(|X_\varepsilon(s)-Z_\varepsilon(s)|)\bigg]ds\cr
	&& \ \ \  \leq \frac{1}{2}  \hat{\mathbb{E}} \bigg[\sup_{0 \leq s \leq u}|X_\varepsilon(s)-Z_\varepsilon(s)|^{2p}\bigg]\cr
	&& \ \ \ +\Upsilon\sqrt{\varepsilon} \int_0^u \lambda(s) \rho_{\eta}\bigg(\hat{\mathbb{E}} \bigg[\sup_{0 \leq r \leq s}|X_\varepsilon(r)-Z_\varepsilon(r)|^{2p}\bigg]\bigg)ds+ \sqrt{\varepsilon} u\Theta,
\end{eqnarray*}

So, we have
\begin{eqnarray*}
	\frac{1}{2} \hat{\mathbb{E}}\bigg[\sup_{0 \leq t \leq u}|X_\varepsilon(t)- Z_\varepsilon(t)|^{2p}\bigg] &\leq& \Upsilon\sqrt{\varepsilon} \int_0^u \lambda(s) \rho_{\eta}\bigg(\hat{\mathbb{E}} \bigg[\sup_{0 \leq r \leq s}|X_\varepsilon(r)-Z_\varepsilon(r)|^{2p}\bigg]\bigg)ds\cr
	&&+ \sqrt{\varepsilon} u\Theta,
\end{eqnarray*}
where $$\Upsilon=((C_{11}+C_{21}+C_{41}+C_{51})\sqrt{\varepsilon}+C_{31}),$$ and $$\Theta= ((C_{12}+C_{22}+C_{42}+C_{52})\sqrt{\varepsilon}+C_{32}).$$

Then by generalization of the Gronwall-Belmman type inequality, we have
\begin{eqnarray*}
	\hat{\mathbb{E}}\bigg[\sup_{0 \leq t \leq u}|X_\varepsilon(t)- Z_\varepsilon(t)|^{2p}\bigg] \leq C {\sqrt{\varepsilon} u\Theta}.
\end{eqnarray*}

Choose suitable $\alpha\in(0,1)$ and $L>0$ such that for every $t\in [0,L\varepsilon^{\frac{1}{2}-\alpha}]\subseteq [0,T],$ we have
\begin{eqnarray*}
	\hat{\mathbb{E}}\bigg[\sup_{0 \leq t \leq L\varepsilon^{\frac{1}{2}-\alpha}}|X_\varepsilon(t)- Z_\varepsilon(t)|^{2p}\bigg]
	\leq Q\varepsilon^{1-\alpha},
\end{eqnarray*}
where $Q=C L\Theta$
is a constant.

Constantly, give any number $\delta_1>0$, we can choose $\varepsilon_1\in(0,\varepsilon_0]$ such that for each $\varepsilon\in(0,\varepsilon_1]$, and for every $t\in [0,L\varepsilon^{\frac{1}{2}-\alpha}]\subseteq [0,T],$ the inequality
\begin{eqnarray*}
	\hat{\mathbb{E}}\bigg[\sup_{0 \leq t \leq L\varepsilon^{\frac{1}{2}-\alpha}}|X_\varepsilon(t)- Z_\varepsilon(t)|^{2p}\bigg]\leq \delta_1.
\end{eqnarray*}
holds. \qed
\para{The Proof of Corollary \ref{thm2}:}
By Lemma \ref{lemcapi} and Theorem \ref{thm1}, for any given number $\delta_2 > 0$, one can find
\begin{eqnarray*}
	\overline{\mathbb{C}}\bigg[\sup_{t\in[0,L\varepsilon^{\frac{1}{2}-\alpha}]}|X_\varepsilon(t)- Z_\varepsilon(t)|>\delta_2\bigg]
	&\leq&  \frac{1}{\delta_2^p} \hat{\mathbb{E}}\bigg[\sup_{0 \leq t \leq L\varepsilon^{\frac{1}{2}-\alpha}}|X_\varepsilon(t)- Z_\varepsilon(t)|^p\bigg]\cr
	&\leq& \frac{Q}{\delta_2^p} \varepsilon^{1-\alpha}.
\end{eqnarray*}
Let $ \varepsilon \rightarrow 0 $ and the required result follows.\qed
\begin{remark}
	{\rm
		When we take $K(t)=0$, (\ref{orginte1}) will reduce to the standard SDEs driven by $G$-Brownian motion and the averaging principle for this equation under condition (H0) has been established (see, e.g. \cite{Han2017}). Different from Han and Liu's work in \cite{Han2017}, we consider the averaging principle for MSDEs driven by $G$-Brownian motion under more general non-Lipschitz condition (H1).}
\end{remark}			
\begin{remark}
	{\rm
		The obtained averaging results in \cite{Guo2018,Mao2019,Xu2011} did not work in sublinear expectation space case, so, in this paper, we obtained the convergence theorem between the solution of averaged multivalued SDEs driven by $G$-Brownian motion and the original one in the sense of mean square (sublinear expectation $\hat{\mathbb{E}}$) and also in capicity which is much more difficulty than the Brownian motion case.}
\end{remark}			
\section{An Example}
Now, we give an example to illustrate our averaging theory.
\begin{example}{\rm
		Consider the following MSDEs driven by $G$-Brownian motion.
		\begin{eqnarray}
		x_\varepsilon(t)&=&X(0)+\varepsilon \int_0^t f(s, x_\varepsilon(s)) ds
		+ \sqrt{\varepsilon}  \int_0^t \sigma(s, x_\varepsilon(s))
		dB_s\cr&&
		+{\varepsilon}\int_{0}^{t}g(s, x_\varepsilon(s))d\langle B\rangle_s- \varepsilon K(t),
		\end{eqnarray}
		where $B_s$ is a $m$-dimensional $G$-Brownian motion and
		\begin{eqnarray*}
		f(s, x) &:=&\sin(s) \sum_{k \geqslant 1} \frac{\sin (k x)}{k^{2}},\cr
		g(s, x)&:=&{\sin(s)}\sum_{k \geq 1} \frac{\sin ^{2}(k x)}{k^{3}}, \cr
		\sigma(s, x) &:=&\sin(s) (1^{-3 / 2} \sin x,  2^{-3 / 2} \sin 2 x, \ldots, m^{-3 / 2} \sin m x ).
			\end{eqnarray*}}
\end{example}

Firstly, we can verify that the functions $f,g,\sigma$ satisfy the condition (H1).
$$
\begin{aligned}|f(t, x)-f(t, y)|  &\leq \sum_{k=1} \frac{|\sin (k x)-\sin (k y)|}{k^{2}}  \\ & \leq C |x-y|\tilde{\kappa}_{1}(x-y|), \end{aligned}
$$
and
$$
\begin{aligned}|\sigma(t, x)-\sigma(t, y)|^{2} &= \sum_{k=1}^{m} \frac{|\sin (k x)-\sin (k y)|^{2}}{k^{3}} \\ & \leq  \sum_{k=1} \frac{|\sin (k x)-\sin (k y)|^{2}}{k^{3}} \\ & \leq 4 \sum_{k=1} \frac{|\sin (k(x-y)) / 2|^{2}}{k^{3}} \\ & \leq C |x-y|^{2} \tilde{\kappa}_{2}(|x-y|), \end{aligned}
$$
and
$$
\begin{aligned} |g(t, x)-g(t, y)|^{2}& \leq C \left(\sum_{k \geq 1}\left|\frac{\sin ^{2}(k x)-\sin ^{2}(k y)}{k^{3}}\right|\right)^{2} \\ & \leq C\left(\sum_{k \geq 1} \frac{|\sin (k x)-\sin (k y)|^{2}}{k^{3}}\right)\\ &
\times \left(\sum_{k \geqslant 1} \frac{|\sin (k x)+\sin (k y)|^{2}}{k^{3}}\right) \\ & \leq C \sum_{k \geq 1} \frac{|\sin (k x)-\sin (k y)|^{2}}{k^{3}} \\ & \leq C |x-y|^{2} \tilde{\kappa}_{2}(|x-y|), \end{aligned}
$$
where
$$
\tilde{\kappa}_{1}(x) :=\left\{\begin{array}{ll}{\log x^{-1},} & {0<x \leq \eta,} \\ {\log \eta^{-1}-1+\frac{\eta}{x},} & {x>\eta,}\end{array}\right.
$$
and $$
\tilde{\kappa}_{2}(x) :=\left\{\begin{array}{ll}{\log x^{-1},}  \quad {0<x \leq \eta,} \\ {\frac{\left(\left(\left(\log \eta^{-1}\right)^{1 / 2}-(1 / 2)\left(\log \eta^{-1}\right)^{-1 / 2}\right) x+(1 / 2)\left(\log \eta^{-1}\right)^{-1 / 2} \eta\right)^{2}}{x^{2}}, \quad x>\eta.}\end{array}\right.
$$
Now, we take $
\kappa_{1}(x) :=C \tilde{\kappa}_{1}(x),\kappa_{2}(x) :=C \tilde{\kappa}_{2}(x)
$. It is easily justified that $\kappa_{1}(x),\kappa_{2}(x)$ satisfy (H1).

Secondly, according the conditions {\rm$(\mathbf{C}_{f}), (\mathbf{C}_{g}), (\mathbf{C}_{\sigma})$}, we can obtain the averaged coefficients.
\begin{eqnarray*}
	\bar{f}(x)&=& \frac{1}{\pi}\int_0^{\pi}{f}(t,x) dt
	=\frac{2}{\pi}\sum_{k \geqslant 1} \frac{\sin (k x)}{k^{2}}\cr
	\bar{\sigma}(x)&=& \frac{1}{\pi}\int_0^{\pi}{\sigma}(t,x) dt=\frac{2}{\pi}{\left(1^{-3 / 2} \sin x, 2^{-3 / 2} \sin 2 x, \ldots, m^{-3 / 2} \sin m x \right)},\cr
	\bar{g}(x)&=& \frac{1}{\pi}\int_0^{\pi}{g}(t,x) dt
	=\frac{2}{\pi} \left(\sum_{k \geq 1} \frac{\sin ^{2}(k x)}{k^{3}}\right).
\end{eqnarray*}

Thus, according to Theorem \ref{thm1}, we can obtain the following convergence result.
\begin{theorem}
	For a given arbitrarily small number $ \delta_1 > 0,$ there exist $ L > 0 $, $\alpha \in (0,1) $, $ \varepsilon_1 \in (0,\varepsilon_0] $, such that for any $ \varepsilon \in (0,\varepsilon_1] $ , each $ t \in
	[0,L\varepsilon^{\frac{1}{2}-\alpha}] $,
	\begin{eqnarray*}
		\hat{\mathbb{E}}\bigg[\sup_{t\in[0,L\varepsilon^{\frac{1}{2}-\alpha}]}|x_\varepsilon(t)- \bar{x}_\varepsilon(t)|^{2p}\bigg]\leq \delta_1.
	\end{eqnarray*}
	where
	\begin{eqnarray*}
		\bar{x}_\varepsilon(t)&=&X(0)+\varepsilon \int_0^t \bar{f} (\bar{x}_\varepsilon(s)) ds
		+ \sqrt{\varepsilon}  \int_0^t \bar{\sigma}  (\bar{x}_\varepsilon(s))
		dB_s\cr
		&&+{\varepsilon}\int_{0}^{t}\bar{g}(\bar{x}_\varepsilon(s))d\langle B\rangle_s- \varepsilon  \bar{K}(t),
	\end{eqnarray*}
	and
	$$
	\bar{f}(x)=\frac{2}{\pi}\sum_{k \geqslant 1} \frac{\sin (k x)}{k^{2}}, \quad \bar{g}(x)=\frac{2}{\pi} \left(\sum_{k \geq 1} \frac{\sin ^{2}(k x)}{k^{3}}\right),$$
	and
	$$\bar{\sigma}(x)=\frac{2}{\pi}{\left(1^{-3 / 2} \sin x, 2^{-3 / 2} \sin 2 x, \ldots, m^{-3 / 2} \sin m x \right)}.$$
\end{theorem}

\section*{Acknowledgments}
This work was partially supported by the National Natural Science Foundation of
China (NSF) under Grant No. 12172285, Guangdong Basic and Applied Basic Research
Foundation under Grant No. 2214050001158 and the Fundamental Research Funds for the
215
Central Universities,China.

\end{document}